\documentclass{amsart}
\usepackage{amsfonts}
\usepackage{amscd}
\usepackage{amsmath, amssymb}

\theoremstyle{plain} \numberwithin{equation}{section}
\newtheorem{theorem}{Theorem}[section]

\newtheorem{lemma}[theorem]{Lemma}
\newtheorem{proposition}[theorem]{Proposition}
\theoremstyle{definition}

\newtheorem{remark}[theorem]{Remark}

\numberwithin{equation}{section}
\newcommand{\1}{\bf 1}
\newcommand{\la}{\langle}
\newcommand{\ra}{\rangle}
\newcommand{\cA}{\mathcal {A}}

\newcommand{\cF}{\mathcal{F}}
\newcommand{\cG}{\mathcal{G}}
\newcommand{\cO}{\mathcal{O}}

\newcommand{\cK}{\mathcal{K}}
\newcommand{\cT}{\mathcal{T}}
\newcommand{\cW}{\mathcal{W}}
\newcommand{\gO}{\Omega}
\newcommand{\cX}{\mathcal{X}}
\newcommand{\ad}{\mathrm{ad}}
\newcommand{\bs}{\backslash}
\newcommand{\fa}{\mathfrak{a}}

\newcommand{\bc}{\bf c}
\newcommand{\fg}{\mathfrak{g}}

\newcommand{\fk}{\mathfrak{k}}
\newcommand{\fm}{\mathfrak{m}}
\newcommand{\fn}{\mathfrak{n}}
\newcommand{\fp}{\mathfrak{p}}

\newcommand{\R}{\mathbb{R}}
\newcommand{\C}{\mathbb{C}}

\newcommand{\Ad}{\mathop{\rm Ad} }
\newcommand{\Sl}{\mathop{\rm Sl} }
\newcommand{\SO}{\mathop{\rm SO} }
\newcommand{\tr}{\mathop{\rm tr} }
\newcommand{\im}{\mathop{\rm im} }

\newcommand{\RE}{\operatorname{Re}}
\newcommand{\IM}{\operatorname{Im}}

\newcommand{\oline}{\overline}

\theoremstyle{plain}

\begin{document}
\title[Heat kernel transform]{The image of the heat kernel
transform on Riemannian symmetric spaces
of the noncompact type}
\author{Bernhard Kr\"otz, Gestur \'Olafsson and Robert J. Stanton}
\address{RIMS, Kyoto University, Kyoto 606-8502, Japan}
\email{kroetz@kurims.kyoto-u.ac.jp}
\address{Department of Mathematics, Louisiana State
University, Baton Rouge,
LA 70803} 
\email{olafsson@math.lsu.edu}
\address{Department of Mathematics, Ohio State 
University, 231 West 18th Avenue, 
Columbus OH 43210-1174}
\email{stanton@math.ohio-state.edu}
\thanks{BK was supported by the RiP-program in Oberwolfach and NSF grant 
DMS-0097314}
\thanks{G\'O was supported by the RiP-program in Oberwolfach,
NSF grants DMS-0139783 and DMS-0402068}
\thanks{RJS was supported in part by NSF grant DMS-0301133} 
\begin{abstract} 
The heat kernel transform
on $G/K$, a Riemannian symmetric space of noncompact type, maps an
$L^2$-function on $G/K$ to a holomorphic function
on the complex crown. In this article we determine
the image of this transform on $L^2$.
\end{abstract}
\maketitle
\section{Introduction}
\noindent
The heat equation, and the associated heat kernel
transform (also known as the Bargmann-Segal transform), is a natural counterpart of a Riemannian metric. For the homogeneous spaces $\mathbb{R}^n$ or a compact symmetric space endowed with the usual metric, many  properties of the heat kernel transform
are known \cite{Barg,H94,S99}. The relevant property for this note is the fundamental observation made by Bargmann and Segal that the image of the heat kernel transform for $\mathbb{R}^n$ consists of functions with holomorphic extension to $\mathbb{C}^n$. Let $X$ denote one of the aforementioned Riemannian manifolds ($\mathbb{R}^n$ or a compact symmetric space) and $k_t$ the heat kernel. The isometry group of $X$ acts transitively, consequently $X$ is diffeomorphic and isometric to $G/K$, $K$ the isotropy group of a basepoint. One knows that the Riemannian manifold $X$ has a natural complexification, denoted here by $\Xi$, which is $G$-diffeomorphic to the tangent bundle $TX$, and on which $k_t$ has a holomorphic extension, say $k_t^\sim$. For a function $f$ on $X$, $k_t*f$, the
solution of the heat equation with initial data $f$, has holomorphic extension to $\Xi$ given by $k_t^\sim*f$. For $f$ in the Hilbert space $L^2(X)$ one obtains thus
a $G$-equivariant linear map,
the \textit{heat kernel transform}, $H_t : L^2(X)\to \cO (\Xi)$.
For these spaces $X$ it is known that the image of $H_t$ is a
weighted Bergman space on $\Xi$, and that $H_t$ is a $G$-equivariant unitary map between these Hilbert spaces. 

\par That this might not be a general 
phenomenon was observed first in \cite{KTX04} for the Heisenberg group with a left invariant Laplace operator. For
this Riemannian space the image of the heat kernel transform is no
longer one weighted Bergman space but a sum
of two, each of which corresponds to a nonpositive weight function 
with strong oscillatory behavior. Surprisingly, neither the Heisenberg group nor the other examples suggest the form of the result for 
symmetric spaces of noncompact type.

\par One problem to be resolved had been to determine the \textit{natural}
complexification of $X=G/K$. Introduced in \cite{AG90} for this purpose, $\Xi$, or as it was later  called the 
{\it complex crown} of $G/K$, 
is a $G$-invariant domain in $X_\C =G_\C/K_\C$ containing $X$, and biholomorphic to an open complex subdomain in $TX$ endowed with the adapted complex structure.  The $G$-orbit structure of $\Xi$ is nicely given by  $\Xi =G\exp (i\gO)\cdot x_o$, where the basepoint $x_o$ is ${eK_\C} \in X_\C$, $X$ is identified with the $G$-orbit through $x_0$, and $\gO$ is a specific polyhedral subset of a maximal flat subspace $\mathfrak{a}\subset T_{x_o}X$. 

In \cite{KS, KS-II} it was shown that $\Xi$ is a maximal domain of holomorphy in $X_\C$ for eigenfunctions of $\mathbb{D}( X)$, the algebra of invariant differential operators on $X$. Consequently, $\Xi$ provides a canonical domain on which to investigate the image of the heat kernel transform on 
$X$. Such a study was initiated in \cite{KS-II} where a
$G$-equivariant holomorphic heat kernel transform for $X$ was
introduced. There it was shown that
the heat kernel $k_t$ extends to a holomorphic
function $k_t^\sim$ on $\Xi$, and that, similarly, one has a $G$-equivariant
heat kernel transform
$H_t : L^2(X) \to \cO (\Xi )$ 
$$H_t (f)(z)= \int_{X} k_t^\sim (g^{-1}z)f(gK)
d\mu_X (gK).$$

The purpose of this note is to give an explicit determination of the image of the heat kernel transform for the noncompact Riemannian symmetric spaces. The description of the image is unlike the previously known cases. That a different type of description is necessary we note in Remark 3.1 (and prove in \S 4) for symmetric spaces with $G$ complex as we show that the image can not be
a weighted Bergman space, i.e., there is no measure of the form $\mu_t = w_t(z) dz$ on $\Xi$ with $w_t$ $G$-invariant, and with 
the norm on $\im H_t$ given by $\int_{\Xi} |F(z)|^2\, d\mu_t(z)$. 

\par Some of the techniques we use were developed in \cite{KS} and \cite{KS-II}. The main new tool, however, is the {\it $G$-orbital integral},
first appearing in \cite{Faraut} and further developed in \cite{KO04},
$$\cO_h(iY) = \int_G h(g\exp ( \frac{i}{2}Y)\cdot x_o)\, dg,$$
$h$ a function on $\Xi$ suitably decreasing at the boundary
and $Y\in 2\gO$. In Theorem \ref{2.3} we show that
there exists a pseudo-differential shift operator
$D$ such that for $f\in L^2(X)$ the function 
$D\cO_{|H_tf|^2}$, initially defined on $2\Omega$, has a natural holomorphic 
extension 
to $\fa_\C$, and that there is an equality of norms 

$$\|f\|^2=\int_{\fa} (D\cO_{|H_tf|^2}) (iY)\, w_t(Y) dY ,$$
the weight function $w_t$ being given by
$$w_t(Y )=\frac{1}{|\cW |}\cdot
\frac{e^{2t|\rho|^2}}{(2\pi t)^{n/2}}\cdot
e^{\frac{-|Y|^2}{2t}}\qquad (Y\in\mathfrak{a}),$$
$\cW$ the usual Weyl group and $n=\mathrm{rank}\, X$. 
To describe $\im H_t$, in \S2 we introduce the space $\cG(\Xi)$ consisting of 
holomorphic functions $F$ on $\Xi$ satisfying the following two  properties:
\begin{itemize}
\item $f:= F|_X \in L^2(X)$
\item $\int_\cX |\widehat f(b,\lambda)|^2\ \psi_\lambda(iY)
\ d\mu(b,\lambda)<\infty$ for all $Y\in \fa$ and $\psi_\lambda(iY)=\sum_{w\in \cW} e^{\lambda(iwY)}$
the Weyl-symmetrized exponential.  
\end{itemize}
Then we show that 
$\im H_t$ is the space
of those functions $F\in \cG(\Xi)$ 
with 
$\int_{\mathfrak{a}} (D\cO_{|F|^2})(iY) w_t(Y)\, dY<\infty$. 

\par We wish to recognize the effort expended by the referee whose questions and insistence on details have led to an improved paper.

\par {\it Special Acknowledgement}: Section 4 and 5 of this paper were 
written while the first named author was visiting 
the R.I.M.S., Kyoto. 
BK would like to express his sincere gratitude to his host 
Toshiyuki Kobayashi for giving him very good advice on style and 
contents of this article.

\section{Riemannian symmetric spaces}  \label{s-one}
\noindent
We recall the necessary background on the geometry and analysis on $\Xi$, the {\it complex crown} of $X$. 

\subsection{The complex crown}\label{ss-one}

Let $G$ be a connected semisimple Lie group. We will 
assume that $G$ is contained in its universal 
complexification $G_\C$ and that $G_\C$ is simply connected. 
Let $K$ be a maximal compact subgroup of $G$ and form the 
homogeneous space $X=G/K$. The complexification $X_\C =G_\C/K_\C$ contains $X$ 
as a totally real submanifold. 

\par Write $\fg$ (resp. $\fk$) for the Lie algebras of $G$
(resp. $K$).
We denote the Cartan-Killing form on $\fg$
by $\la\cdot, \cdot\ra$. Then $\fp :=\fk^\bot$ is a complementary
subspace on which $\la\cdot, \cdot\ra$ defines an inner product and hence a Riemannian symmetric structure on $X=G/K$.
Fix a maximal abelian subspace $\fa\subset \fp$ and set
$$\Omega=\left\{ Y\in \fa: {\rm Spec} (\ad Y)\subset \left]-
\frac{\pi}{2}, \frac{\pi}{ 2}\right[\right\}.$$
$\Omega$ is open, convex, and invariant
under the  Weyl group $\cW=N_K(\fa)/M$, $M=Z_K(\fa)$.

\par Let $x_o=eK_\C$ be the base point in $X_\C$. The 
{\it complex crown of $X$} is defined by 
$$\Xi=G\exp(i\Omega)\cdot x_o.$$
From several sources (for example \cite{BHH} or \cite{KS-II} and references therein) 
one knows that $\Xi$ is a $G$-invariant Stein 
domain in $X_\C$ which contains $X$ in its ``middle'' as 
a totally real submanifold. 
The definition of $\Xi$ is 
independent of the choice of $\fa$ (all maximal
abelian subspaces in $\fp$ are conjugate under
$K$), hence $\Xi$ is canonically defined by $X$.

\par Denote
by $\Sigma\subset \fa^*$ the set of restricted roots. Then

$$\Omega=\left\{ Y\in \fa: (\forall \alpha\in \Sigma)\ |\alpha(Y)|<
\frac{\pi}{2}\right\}. $$
In particular $\oline \Omega$ is a compact  polyhedron. 
\par To each $\alpha\in \Sigma$ we associate the root space 
$\fg^\alpha=\{ Y\in \fg: (\forall H\in \fa)\  [H,Y]=\alpha(H) Y\}$. 
If $\fm$ denotes the Lie algebra of $M$, one has the root space 
decomposition 
$$\fg=\fa\oplus\fm\oplus\bigoplus_{\alpha\in \Sigma}\fg^\alpha.$$
Fix a positive system $\Sigma^+\subset \Sigma$ and set 
$\fn:=\bigoplus_{\alpha\in \Sigma^+} \fg^\alpha$. If $A:=\exp (\fa)$ and
$N:=\exp (\fn) $ are the (closed) analytic subgroups
corresponding to $\fa$ and $\fn$,
and $A_\C:=\exp (\fa_\C)$ and $N_\C:=\exp (\fn_\C)$ are
the corresponding complex groups, then one has the inclusion
$$\Xi\subset N_\C A_\C \cdot x_o$$
(cf.\  \cite{KS}, \cite{Huck}, \cite{Mat}).
As $\Xi$ is contractible one can define
an $N$-invariant holomorphic map
$$a: \Xi \to A_\C$$
with $a(x_o)=\1$ and having $z\in N_\C a(z)\cdot x_o$
for all $z\in \Xi$. The image of
$\IM \log \circ a$ is described by the complex convexity
theorem (cf.\ \cite{GK}, \cite{KO}):

\begin{equation}\label{eq=cc}
\IM \log a(G\exp (iY)\cdot x_o)={\rm conv}(\cW\cdot Y).
\end{equation}

\subsection{Spherical functions}\label{ss-two}

Spherical functions on $X$ have a
holomorphic extension to the $G$-domain $\Xi$ (cf.\ \cite{KS}). Moreover, $\Xi$
is a maximal $G$-domain in $X_\C$ with this property (see 
Theorem \ref{lasto-one} below). Thus, it seems reasonable 
to define 
spherical functions ab initio as holomorphic functions 
on $\Xi$. Detailed proofs of the material below can be found 
in \cite{KS}.

\par Set $\rho=\frac{1}{ 2}\sum_{\alpha\in \Sigma^+} (\dim \fg^\alpha)
\alpha\in \fa^*$. For $\lambda\in \fa_\C^*$ we set
$a(\cdot)^\lambda :=e^{\lambda \log a(\cdot)}$. 
Notice that  $a^\lambda$ is a holomorphic function on $\Xi$.

\par The {\it spherical function with parameter} $\lambda\in \fa_\C^*$ is the 
holomorphic
function on $\Xi$ defined by

$$\varphi_\lambda(z)=\int_K a(kz)^{\rho+\lambda}\ dk\qquad (z\in \Xi),$$
$dk$, as usual, a normalized Haar measure on the compact 
group $K$. 

\par Later we shall restrict our attention to spherical functions with 
imaginary 
parameters, i.e. $\lambda\in i\fa^*$, as only these appear in the 
Plancherel decomposition of 
$L^2(X)$.

\par While $\Xi$ is canonical as a $G$-domain in $X_\C$ for spherical functions, the situation is slightly different if one considers 
instead the action of $K_\C$ on $X_\C$. 
Write $\widehat X_{\C,2\Omega}$ for the smallest $K_\C$-invariant 
domain in $X_\C$ which contains $A\exp(i2\Omega)\cdot x_o$ (for existence see 
\cite{K}). Then a 
spherical function $\varphi_\lambda$  can be considered equally
as a $K_\C$-invariant holomorphic function on $\widehat X_{\C,2\Omega}$. 
To avoid unnecessary notation, we write $\varphi_\lambda$ also when 
viewed as a function on $\widehat X_{\C, 2\Omega}$. 
For $u\in A$ and $Y\in \Omega$ we have then 

\begin{equation}\label{eq=se} 
|\varphi_\lambda(u\exp(i2Y)\cdot x_o)|\leq 
\varphi_\lambda(\exp(i2Y)\cdot x_o)=\int_K 
\left|a(k\exp(iY)\cdot x_o)^{\rho+\lambda}\right|^2\ dk . 
\end{equation}

\par For $\lambda\in i\fa^*$ we set 
\begin{equation}\label{eq=psi} 
\psi_\lambda(Z)=\sum_{w\in \cW} e^{\lambda(wZ)} \qquad (Z\in \fa_\C).
\end{equation}
Notice that $\psi_\lambda$ is a $\cW$-invariant function that is 
positive on $i\fa$ and of exponential growth there. 

\begin{lemma} For each $Y\in \gO$ there exists a constant
$C_Y>0$ such that
\begin{equation}\label{eq=esa}
\varphi_\lambda(\exp(i2Y)\cdot x_o)\leq C_Y\cdot  \psi_\lambda(i2Y)
\end{equation}
for all $\lambda\in i\fa^*$.
\end{lemma}

\begin{proof}This follows from equations (\ref{eq=cc}) and (\ref{eq=se}).
\end{proof}

\begin{remark} $C_Y$ is locally bounded on $\Omega$.
\end{remark}
 
\subsection{Harmonic analysis on $\Xi$}  \label{ss=three}
We conclude this section with some aspects of harmonic analysis
on $\Xi$. We begin with the Fourier transform 
and Plancherel theorem for $X$ following Helgason and Harish-Chandra. 

\par Let $B=M\bs K$ and $db$ a normalized Haar measure on $B$. Set $\cX=B\times i\fa^*$ and define a measure $\mu$ on 
$\cX$ by $d\mu(b,\lambda)= db\otimes \frac{d\lambda}{|\bc(\lambda)|^2}$ where
$d\lambda$ is the measure on $i\fa^*$ as normalized by Harish-Chandra, $\bc(\lambda)$ is the Harish-Chandra
$\bc$-function. For an integrable
function $f$ on $X$ we define its Fourier transform to be 
the function $\widehat f$ on $\cX$ defined by 

$$\widehat f(b, \lambda)=\int_X f(x)\ a(bx)^{\rho-\lambda} \ dx,\qquad 
(b,\lambda)\in \cX,$$
where $dx$ is a $G$-invariant measure on $X$ normalized so that 
$f\mapsto \widehat f$ extends to an isometry $L^2(X)\to L^2(\cX, \mu)$.
The Fourier inversion theorem in this formulation becomes
$$f(x)=\int_\cX \widehat f(b,\lambda)\ a(bx)^{\rho+\lambda} \ d\mu(b,\lambda).$$

\par
Let $\cG(\Xi)$ be the space of holomorphic functions
$F$ on $\Xi$ such that the following two  properties hold:
\begin{itemize}
\item $f:= F|_X \in L^2(X)$
\item $\int_\cX |\widehat f(b,\lambda)|^2\ \psi_\lambda(iY)
\ d\mu(b,\lambda)<\infty$ for all $Y\in \fa$. 
\end{itemize}

\begin{lemma}[Gutzmer's identity d'apr\`es Faraut]\label{le-gi} Let
$F\in \cG(\Xi)$. Then
\begin{equation}\label{eq=Gutzmer}
\int_G |F(g\exp(iY)\cdot x_o)|^2 \ dg =\int_{\mathcal X} |\widehat f(b,\lambda)|^2 
\ \varphi_\lambda(\exp(i2Y)\cdot x_o) \ d\mu(b,\lambda)
\end{equation}   
for all $Y\in \Omega$. 
\end{lemma}
\begin{proof} This is the Gutzmer identity of  \cite{Faraut}, Th.\ 1. 
\end{proof}

\par For a sufficiently decreasing function $h$ on
$\Xi$ we consider its {\it $G$-orbital integral}

$$\cO_{h}(iY)=\int_G h(g\exp(\frac{i}{2}Y)\cdot x_o) \ dg\qquad (Y\in 2\Omega).$$

\noindent In particular, if $F\in \cG(\Xi)$  one has from (\ref{eq=Gutzmer}) that $\cO_{|F|^2}(iY)$ is finite for all
$Y\in 2\Omega$. We consider in $\fa_\C$ the abelian tube domain  $\cT(2\Omega)=\fa +i2\Omega$. Then for
$F\in \cG(\Xi)$ it follows from Lemma \ref{le-gi} and the holomorphic extension of spherical functions \cite{KS} that
$\cO_{|F|^2}$ admits a natural extension to
a holomorphic function on $\cT(2\Omega)$, namely

\begin{equation} \label{eq=exten} \cO_{|F|^2}(Z)=
\int_{\mathcal X} |\widehat f(b,\lambda)|^2
\ \varphi_\lambda(\exp(Z)\cdot x_o) \ d\mu(b,\lambda)\qquad (Z\in \cT(2\Omega)).
\end{equation}

\section {The heat kernel transform}\label{s-two}
\noindent
In this section we shall
determine the image of the heat kernel
transform on $L^2(X)$ and its holomorphic extension
to $\Xi$.
We start with a brief review of the results from \cite{KS-II}.

\subsection{Definition and basic properties}\label{ss-two-one}

In the following, $t$ denotes a positive number. According to Gangolli 
\cite{Gang}
the heat kernel $k_t$, suitably normalized,  on $X$ has the spectral resolution 

\begin{equation}\label{eq=heat} k_t(x)= \int_{i\fa^*} e^{-t(|\lambda|^2 
+|\rho|^2)}
\varphi_\lambda(x) \ \frac{d\lambda}{ |\bc(\lambda)|^2}\qquad (x\in X).
\end{equation}

\par If $f$ is an analytic function on $X$ which admits a holomorphic extension 
to $\Xi$, then we write  $f^\sim$ for the holomorphically extended function.   
It is proved in \cite{KS-II} (but also immediate from (\ref{eq=heat}) and the 
results collected in Subsection \ref{ss-two}) that $k_t$ admits an analytic 
continuation 
$k_t^\sim$
given by 

$$k_t^\sim(z)= \int_{i\fa^*} e^{-t(|\lambda|^2 +|\rho|^2)}
\varphi_\lambda(z) \ \frac{d\lambda}{ |\bc(\lambda)|^2}\qquad (z\in \Xi).$$
Our concern here is the heat kernel transform, viz. the operator given by convolution of $k_t$ on $L^2(X)$: 

$$(k_t*f)(x)=\int_X k_t(g^{-1}x) f(gK) \ d(gK)\qquad (x\in X).$$ 
It was proved in \cite{KS-II} that $k_t*f$ admits a holomorphic extension 
given  by 
$$(k_t*f)^\sim (z)=\int_X k_t^\sim (g^{-1}z) f(gK) \ d(gK)\qquad (z\in \Xi).$$ 
Thus one gets a linear map 

$$H_t: L^2(X)\to \cO(\Xi), \ \ f\mapsto (k_t*f)^\sim 
$$
which is refered to as the {\it heat kernel transform with 
parameter $t>0$}. In what follows
we consider $\cO(\Xi)$ as a topological vector space endowed 
with the usual Fr\'echet topology of compact convergence. 

In \cite{KS-II} the following properties of the heat kernel transform are proved:

\begin{itemize}
\item $H_t$ is continuous;
\item $H_t$ is $G$-equivariant;
\item $H_t$ is injective.
\end{itemize}

\noindent Thus $\im H_t$, when endowed with the
Hilbert topology from $L^2(X)$, becomes a $G$-invariant Hilbert space
of holomorphic functions on $\Xi$. As such it admits  a reproducing kernel 
which was spectrally characterized in (\cite{KS-II}) although the image of $H_t$ was not
described in \cite{KS-II}. 

\begin{remark}\label{ne}  Contrary to expectations, we suspect that $\im H_t$ is not a 
weighted Bergman space,
i.e. there does \textit{not} exist a measure $\mu_t$ on $\Xi$, 
absolutely continuous with respect to the standard Lebesgue measure,
such that the norm on $\im H_t$ is given by
$$\int_\Xi |F(z)|^2 \ d\mu_t(z) \qquad (F\in \im H_t).$$
In Section \ref{NE} we will present an argument that shows for
complex groups that there does not exist a $G$-invariant $\mu_t$.  
\end{remark}

\subsection{The image of the heat kernel transform}\label{ss-two-two}
\par Let $f\in L^2(X)$ and set $F=H_t(f)$. In \cite{KS-II} it was shown that $F\in \cO (\Xi)$. As $\psi_\lambda$ is only of exponential growth it is easily seen that $F\in \cG(\Xi)$. In particular, 
we can form the orbital integral $\cO_{|F|^2}$ (cf.\ (\ref{eq=exten})) and  
for all $Z\in \cT(2\Omega)$ one has
from Lemma 2.3
\begin{eqnarray}\label{eq=orbit} \cO_{|F|^2} (Z) &=& \int_\cX |\widehat 
F(b,\lambda)|^2\ 
\varphi_\lambda(\exp (Z)\cdot x_o)
\ d\mu(b,\lambda) \\
&=& \int_\cX |\widehat f(b,\lambda)|^2 \ e^{-2t(|\lambda|^2 +|\rho|^2)}
\ \varphi_\lambda(\exp(Z)\cdot x_o)
\ d\mu(b,\lambda) \\
&=& \int_{i\fa^*} g(\lambda) \ \varphi_\lambda(\exp(Z)\cdot x_o)
\ \frac{d\lambda}{ |\bc(\lambda)|^2}\ ,
\end{eqnarray}
where the function $g$ is given  by 
\begin{equation}\label{eq=g}
g(\lambda)=e^{-2t(|\lambda|^2 +|\rho|^2)} \int_B |\widehat  f(b,\lambda)|^2 \ db. 
\end{equation}

\medskip Next we define 
$\cF(\cT(\Omega))$ as the space 
of holomorphic functions $h$ on $\cT(2\Omega)$ such that 
$$h(Z)=\int_{i\fa^*} g(\lambda)\  \varphi_\lambda(\exp(Z)\cdot x_o)\  
\frac{d\lambda}{ 
|\bc(\lambda)|^2}
\qquad (Z\in \cT(2\Omega))$$
for an integrable  ${\mathcal W}$-invariant function $g$ on $i\fa^*$ with 
$\int_{i\fa^*} |g(\lambda)|\  \psi_\lambda(iY)\  \frac{d\lambda}{ 
|\bc(\lambda)|^2}<\infty$
for all $Y\in \fa$. 
Notice that $g$ is uniquely determined by 
$h$. Thus we can define a linear operator 
$$D: \cF(\cT(\Omega))\to \cO(\fa_\C), \ \ (Dh)(Z)=
\int_{i\fa^*} g(\lambda)\  \psi_\lambda(Z) \  \frac{d\lambda}{ 
|\bc(\lambda)|^2},$$
where $\psi_\lambda$ is defined in \ref{eq=psi}.
Observe that for $f\in L^2(X)$ and $F=H_t(f)$ we have 
$\cO_{|F|^2}\in \cF(\cT(2\Omega))$ by (\ref{eq=orbit}) and (\ref{eq=g}),  
and so $D\cO_{|F|^2}\in \cO(\fa_\C)$. 

\begin{remark} In general, the operator $D$ is a pseudo-differential shift operator, but if the multiplicities of all the restricted roots are even, it is a differential shift operator. This is most easily seen by relating the operator $D$ 
to the Abel transform (Harish-Chandra's map $f\to F_f$). 
So let $h\in \cF(\cT(2\Omega))$. Then $h\circ \log|_A$ is a
$\cW$-invariant smooth function on $A$. By the $C^\infty$ Chevalley theorem (cf.\ \cite{Dadok}) it has a natural extension to a $K$-invariant smooth function on $X$, say $H\in C^\infty(X)^K$.
From the growth restriction on $g$
it follows
that $H$ belongs to the Harish-Chandra Schwartz space on $X$. 
Hence the  Abel transform 
\begin{equation} \label{eq=Abel}
(\cA H)(u)=u^\rho\int_N H(un)\ dn \qquad (u\in A)
\end{equation}
of $H$ is  well defined, and $\cA H$ is a $\cW$-invariant 
smooth function on $A$. (We remark that (\ref{eq=Abel}) cannot be extended to all of 
$u\in \exp(\cT(2\Omega))$.) 

Using the well known identity of the spherical transform of $H$ and the Fourier transform of $\cA H$ (for example cf.\ \cite{Sawyer}) one obtains 
$$(\cA H) (u)=\int_{i\fa^*} g(\lambda)\  \psi_\lambda(\log u) \ 
d\lambda \ .$$  

This shows that $D$ is given by $\cA$ composed 
with the Fourier multiplier operator corresponding to $\frac{1}{|\bc(\lambda)|^2}$. The remark follows from the well known properties of $\cA$ and the structure of $\frac{1}{|\bc(\lambda)|^2}$.

\end{remark}

To formulate the main result of the paper, with $n=\dim \fa$, we define a weight function on $\fa$ by 

$$w_t(Y)=\frac{1}{ |\cW|} \cdot\frac{e^{2t|\rho|^2}}{ (2\pi t)^{\frac{n}{ 2}}}
\cdot  e^{-\frac{|Y|^2}{2t}}
\qquad (Y\in \fa).$$

\begin{theorem}\label{2.3} For $f\in L^2(X)$ and $F=H_t(f)$ the following norm 
identity holds
\begin{equation}\label {eq=norm}\|f\|^2 =\int_{\fa} (D\cO_{|F|^2})(iY)\  w_t(Y) 
\ dY.
\end{equation}
Moreover, one has 
\begin{equation}\label{eq=char}\im H_t =\{ F\in \cG(\Xi): \int_{\fa} 
(D\cO_{|F|^2})(iY)
\  w_t(Y) \ dY<\infty\}.
\end{equation}
\end{theorem}
\begin{proof} As observed earlier, $F=H_t(f)$ is in $\cG(\Xi)$. So consider the norm identity (\ref{eq=norm}). 
By the Plancherel theorem the left hand side of (\ref{eq=norm}) is given by 
\begin{equation} \label{eq=left}
\|f\|^2=\int_\cX |\widehat  f(b,\lambda)|^2 \ d\mu(b,\lambda)
\end{equation}
while, according to (\ref{eq=orbit}), the right hand side of (\ref{eq=norm})  is given by 
\begin{equation}\label{eq=right}
\int_{\fa} (D\cO_{|F|^2})(iY) \ w_t(Y) \ dY
=\int_\fa \int_\cX |\widehat f(b,\lambda)|^2\  e^{-2t(|\lambda|^2 +|\rho|^2)}
\ \psi_\lambda(i2Y) \ w_t(Y)\ d\mu(b,\lambda)\ dY.
\end{equation} 
Comparing (\ref{eq=left}) and (\ref{eq=right}), we see
that (\ref{eq=norm}) will be established 
provided we can show

$$\int_\fa \psi_\lambda(i2Y)\  w_t(Y) dY =e^{2t(|\lambda|^2+|\rho|^2)}$$
for $\psi_\lambda(i2Y)$ the Weyl symmetrized exponential function, and for each $\lambda\in i\fa^*$. 
But, using the definition of $w_t$, we have
\begin{eqnarray*}
\int_\fa \psi_\lambda(i2Y)\  w_t(Y) \ dY &=&
\frac{1}{|\cW|} \cdot  \frac{e^{2t|\rho|^2} }{ (2\pi t)^{\frac{n}{ 2}}}
\int_{\fa} \psi_\lambda(i2Y)\ e^{\frac{-|Y|^2}{ 2t}}\ dY\\
&=& \frac{e^{2t|\rho|^2}}{(2\pi t)^{\frac{n}{ 2}}}
\int_{\fa} e^{\lambda(i2Y)}\ e^{\frac{-|Y|^2}{ 2t}}\ dY\\
&=& \frac{e^{2t|\rho|^2}}{(2\pi )^{\frac{n}{2}}}
\int_{\fa} e^{\lambda(i2\sqrt{t}Y)}\ e^{\frac{-|Y|^2}{ 2}}\ dY\\
&=& e^{2t|\rho|^2} e^{\frac{|2\sqrt{t}\lambda|^2}{ 2}}
=e^{2t(|\lambda|^2 +|\rho|^2)}\ ,
\end{eqnarray*} 

\noindent establishing (\ref{eq=norm}) and also proving that $H_t$ is an injective
isometry into the space described in (\ref{eq=char}).
\par For surjectivity, take an $F\in \cG(\Xi)$.
Then $f= F|_X \in L^2(X)$. We shall use
the Paley-Wiener space from \cite {KS-II}. The main properties of
$PW(G/K)$ needed here are that it is dense in $L^2(X)$ and that
it consists of functions with compactly supported Fourier transforms.
We will assume that $f$ is in $PW(G/K)$. Multiply its Fourier transform
by $e^{t(|\lambda|^2+|\rho|^2)}$, then invert the result.  We obtain a
function $g$ still in $PW(G/K)$. To $g$ we apply the heat transform $H_t$.
Then $H_t(g)$ is in $\cG(\Xi)$ and its restriction to $X$ has the same
Fourier transform as $f$. Since $H_t(g)$ and $F$ are in $\cG(\Xi)$ and their
restrictions to $X$ have the same Fourier transform, we have $H_t(g) = F$.
To conclude the argument we use that $PW(G/K)$ is dense and the identity
of norms in (\ref{eq=norm}). Finally, the characterization (\ref{eq=char})
follows from the observation that only equalities have been used in the above proof.
\end{proof}

\section{$\im H_t$ need not be weighted}\label{NE}

In this appendix we will prove the assertion mentioned 
in Remark \ref{ne}. The amount of detail presented is at the insistence of the referee. 
The material is extracted from 
unpublished notes of the first and last named authors.

\smallskip In many ways, Riemannian symmetric spaces of the non-compact 
type $X=G/K$ resemble the real numbers $\R$, e.g. both $X$ and 
$\R$ are analytic complete Riemannian manifolds with 
Laplace operator having  purely continuous $L^2$-spectrum. There 
is however a basic difference between their natural 
complexifications, $\C=T\R$ and $\Xi$, in that for the fibration 
$\C=T\R \to \R$ the fibers are complete and $T\R=\C$ is a complete manifold; 
whereas, $\Xi=G\times_K \Ad(K)\Omega \subset TX=G\times_K\fp$ resembles 
is a disk bundle over $X$ with neither the fibers nor the total space $\Xi$ complete with respect to the natural (adapted) metric
(\cite{KS-II}, Sect. 4). 
Thus in the abelian setting a better analog of the crown is  a strip domain $S=\{z\in \C: |\IM z|< \gamma\}$, $\gamma>0$. 
For this simple situation it is 
not difficult to see that the image of the 
heat kernel transform $H_t: L^2(\R)\to \cO(\C)$,  when considered 
in $\cO(S)$, is not a weighted Bergman space. 
We include a detailed discussion of the flat case, then
move to the curved situation $\Xi$.

\subsection{Some features of the heat kernel transform on $\R$}

\par Write $\Delta=\frac{d^2}{dx^2}$ for the Laplace operator 
on $\R$ and consider the heat equation 
$$(\Delta -\partial_t) u(x,t)=0 \qquad \left((x,t)\in \R\times \R_{>0}\right).$$ 
The fundamental solution is given by the heat 
kernel 

\begin{eqnarray} k_t(x)&=&
\frac{1}{ \sqrt{2\pi}} \int_\R e^{ixy} e^{-ty^{2}}\ dy \\
\label{heat=R}&=& \frac{1}{\sqrt{4\pi t}} e^{-\frac{x^2}{4t}}
\end{eqnarray}
where $x\in \R, t>0$. 
 
\par If $f$ is an analytic function on $\R$ which extends holomorphically 
to $\C$, then we denote by $f^\sim$ this holomorphic extension. 
From the explicit formula (\ref{heat=R}) one sees that $k_t$ extends holomorphically to $\C$, and 
 
\begin{equation} \label{heat=C}
k_t^\sim(z)=\frac{1}{\sqrt{4\pi t}} \cdot e^{-\frac{z^2}{ 4t}} \qquad (z\in \C).
\end{equation}
 
Fix $f\in L^2(\R)$. The convolution 

\begin{eqnarray*} (f*k_t)(x)&=&\int_\R f(y) k_t(x-y)\ dy \\
&=&\frac{1}{\sqrt{4\pi t}} \int_\R f(y) e^{- \frac{(x-y)^2}{ 4t}} dy 
\end{eqnarray*}
admits a holomorphic continuation to $\C$. Thus we obtain a map 

$$H_t: L^2(\R)\to {\mathcal O}(\C), \ \ f\mapsto (f*k_t)^\sim ,$$
the {\it heat kernel transform}.
We will consider ${\mathcal O}(\C)$ as a Fr\'echet topological 
vector space (topology of compact convergence).
 
It is an elementary task to verify the following properties 
of $H_t$:

\medskip 
\begin {itemize} 
\item $H_t$  is injective.
\item $H_t$  is continuous.
\item $H_t$  is equivariant with respect to the $\R$-action, 
i.e. $H_t(f(\cdot +x))= H_t(f)(\cdot +x)$ for all $x\in\R$ and 
$f\in L^2(\R)$. 
\end{itemize}
\medskip 

From these properties one deduces that ${\mathcal F}_t(\C):=\im H_t$, 
when endowed with the Hilbert inner product from $L^2(\R)$,  
becomes a Hilbert space of holomorphic functions.   
As such it admits a reproducing kernel 
$${\mathcal K}^t: \C\times \C \to \C, \qquad (z,w)\mapsto {\mathcal K}^t(z,w).$$
This is a continuous function which is holomorphic in the first variable
and anti-holomorphic in the second variable. For $w\in \C$ 
we define a function ${\mathcal K}_w^t\in {\mathcal F}_t(\C)$ by 
$${\mathcal K}_w^t(z)={\mathcal K}^t(z,w).$$
If $\la\cdot, \cdot\ra$ denotes the Hilbert inner product 
on ${\mathcal F}_t(\C)$, then ${\mathcal K}_w^t$ describes the point-evaluation 
at $w$, i.e. 

$$(\forall f\in {\mathcal F}_t(\C)) \qquad f(w)=\la f, {\mathcal K}_w^t\ra .$$
Notice that $H_t^{-1}({\mathcal K}_z^t) =\oline {k_t^\sim(z-\cdot)}$.
Thus 

\begin{eqnarray}\label{Kern} {\mathcal K}^t(z,w)&=&\la {\mathcal K}_w^t, {\mathcal K}_z^t\ra =\la H_t^{-1}({\mathcal K}_w^t), 
H_t^{-1} ({\mathcal K}_z^t)\ra_{L^2(\R)}\\ 
\notag &=& \int_\R k_t^\sim (z-x) \oline{k_t^\sim (w-x)} \ dx \\
\notag &=&\int_\R k_t^\sim (z-x) k_t^\sim (\oline w-x) \ dx \\ 
\notag &=&\int_\R k_t^\sim (-x) k_t^\sim (\oline w -z-x) \ dx \\ 
\notag &=&\int_\R k_t^\sim (x) k_t^\sim (\oline w -z-x) \ dx \\ 
\notag &=&(k_t*k_t)^\sim (\oline w-z)\\
\notag &=& k_{2t}^\sim (z-\oline w)=
\frac{1}{ \sqrt{8\pi t}}\cdot  e^{- \frac{(z-\oline w)^2}{8t}}.
\end{eqnarray}

The natural question as to whether ${\mathcal F}_t(\C)$ admits 
a description as a weighted Bergman space was answered
affirmatively
by Bargmann and Segal. 
With 

\begin{equation}\label{weight} 
w_t(y)=\frac{1}{\sqrt{2\pi t}} e^{-\frac{y^2}{ 2t}} \qquad (y\in \R)
\end{equation}
one has 

$${\mathcal F}_t(\C)=\left\{ f\in {\mathcal O}(\C): 
\|f\|^2= \int_\C |f(x+iy)|^2 w_t(y)
\ dxdy<\infty \right\}.$$

This can be easily verified either by direct 
computation or, perhaps preferably, by abstract 
Hilbert space techniques (use that the point evaluations
$\{ {\mathcal K}_z^t: z\in \C\}$ form a dense subspace of ${\mathcal F}_t(\C)$). 

\par Let us pass now to a strip domain $S=\{ z\in \C\mid |\IM z|<\gamma\}$
in the complex plane. 
We denote by ${\rm Res}: \cO(\C)\to \cO(S)$ the restriction map and 
consider the restriction of the heat kernel transform 
${\rm Res}\circ H_t: L^2(\R)\to \cO(S)$. We notice, 
as ${\rm Res}$ is continuous and $S$ is translation invariant, 
that ${\rm Res\circ H_t}$ satisfies 
the bulleted items from before as well. In conclusion, 
$\im {\rm Res}\circ H_t=\cF_t(S)$ is an $\R$-invariant 
Hilbert space of holomorphic functions on $S$. 
However, the nature of  $\cF_t(S)$ is different from 
$\cF_t(\C)$ because 

\begin{proposition} There does not exist a positive 
$\R$-translation invariant measurable weight 
function $w_t: S \to \R_{>0}$
such that 
$$\cF_t(S)=\{ f\in \cO(S): \|f\|^2= \int_S |f(z)|^2 w_t(z)\,  dx\, dy 
<\infty\}; $$
in other words, $\cF_t(S)$ is not a weighted Bergman space 
with respect to a translation invariant weight. 
\end{proposition}

\begin{proof} The reproducing kernel for $\cF_t(S)$ is 
simply the restriction of the kernel of $\cF_t(\C)$  to 
$S\times S$ -- a quick inspection of our derivation   of
the identity $\cK^t(z,w) =
k_{2t}^\sim (z - \overline{w})$ shows that. 
To simplify notation, we write $H_t$ instead of 
${\rm Res}\circ H_t$.  
\par  For each $z\in S$ the function 
$w\mapsto \cK_z^t(w) = \cK^t(w,z)$ is in the Hilbert space
$\im H_t$.  As $k_t*k_t=k_{2t}$, 
we have $\cK_x^t= H_t(k_t(\cdot -x))$ for all $x\in\R$.  
Notice that the collection of $\cK_x^t,\ x\in \R$, spans a dense subspace
of $\cF_t(S)$, hence the existence 
of a translation invariant weight for $\cF_t(S)$ is equivalent to 
\begin{eqnarray}\label{e1}
\langle k_t(\cdot -x), k_t(\cdot -y)\rangle_{L^2(\R)} & = & 
\int_S \cK_x^t(z) \overline {\cK_y^t(z)} \ w_t (z) dz \\ 
\notag & = &  
\int_S \cK_x^t(u+iv) 
\overline {\cK_y^t(u+iv)} \ w_t (v)\, dudv\qquad (x,y\in \R). 
\end{eqnarray}

Using translation invariance we may assume $y=0$. 
Straightforward computation 
transforms the identity (\ref{e1}) into 
\begin{equation}\label{e2}
k_{2t}(x)= \int_{-\gamma}^\gamma
 k_{4t}^\sim (x-2iv) w_t(v) dv \qquad (x\in \R)\, .
\end{equation}
As $k_{4t}^\sim (x-2iv) = {\rm const} 
\cdot k_{4t}(x) e^{ixv/4t} e^{v^2/4t}$, moving $k_{4t}(x)$ to the left side of (\ref{e2})
gives 

\begin{eqnarray*} e^{-{x^2/16t}} &=& {\rm const} \cdot 
\int_{-\gamma}^\gamma  e^{ixv/4t} e^{v^2/4t} w_t(v)\, dv \\ 
e^{-{ty^2}}&=& {\rm const} \cdot 
\int_{-\gamma}^{\gamma}  e^{iyv} e^{v^2/4t} w_t(v)\, dv \qquad 
(y\in \R),
\end{eqnarray*}
and  we obtain a contradiction to the 
Fourier transform of the function  $e^{-ty^2}$.
\end{proof}

\subsection{From flat to curved}

Now moving from $\R$ to $X$ 
we will show that $\cF_t(\Xi)\:=\im H_t $ is NOT a weighted Bergman space.

\par At first we make no restriction on $G$, such as 
$G$ is complex. We recall from \cite{KS-II}, Th. 6.4 the formula for the reproducing 
kernel of $\cF_t(\Xi)$, 

\begin{eqnarray*} \cK^t(z,w) &=& 
\int_X k_t^\sim(g^{-1}z) \oline{k_t^\sim(g^{-1}w)} \ d(gK)\\
&=& \int_X k_t^\sim(g^{-1}z) k_t^\sim(g^{-1}\oline w) \ d(gK) \qquad (z,w\in \Xi).
\end{eqnarray*}
When $\oline w^{-1}z$ is in $K_\C \Xi$ we can 
use integration by parts and $k_t*k_s=k_{t+s}$ to arrive at 
 
\begin{equation}\label{3.2} \cK^t(z,w)=k_{2t}^\sim (\oline w^{-1}z) \qquad (\forall z,w\in \Xi
\ \hbox{with }  \ \oline w^{-1}z\in K_\C \Xi).
\end{equation}
 Then (\ref{3.2}) is the curved analogue of the 
formula (\ref{Kern}). 

\par Suppose for a moment that 
$\cF_t(\Xi)$ were a weighted 
Bergman space with respect to an absolutely continuous 
$G$-invariant measure $\mu_t$, i.e.

$$\cF_t(\Xi)=\{f\in \cO(\Xi) : \|f\|^2=
\int_\Xi |f(z)|^2 \, d \mu_t(z)\ dz<\infty\} . $$

\par As $X\subseteq \Xi$ is a totally real submanifold, it follows that 
the real point evaluation $K_x^t$, $x\in X$, form a 
dense subspace of $\cF_t(\Xi)$. Thus the 
measure $\mu_t$ would be uniquely determined by the values 

\begin{equation}\label{3.3}  \la  \cK_x^t, \cK_y^t\ra =\int_\Xi \cK_x^t(z) \oline {\cK_y^t(z)} \, d\mu_t(z)
\end{equation}
for $x, y\in X$. Using $G$-invariance we may actually assume 
that $y=x_o$ and $x=a\cdot x_o$ for some $a\in A$. 
Apply (\ref{3.2}) and transform the left hand side of (\ref{3.3}) 
to  
 
\begin{equation} \label{3.4} 
\la \cK_x^t, \cK_y^t\ra= \cK^t(a\cdot x_o, x_o)=k_{2t}(a^{-1})=k_{2t}(a) \qquad (a\in A).
\end{equation}

We will write $\Omega^+$ for the intersection of $\Omega$ with 
a Weyl chamber and recall from \cite{KS-II}, Cor. 4.2,  the fact that the map 
$$G/M\times \Omega^+\to \Xi, \ \ (gM, Y)\mapsto g\exp(iY)\cdot x_o$$
is a $G$-equivariant diffeomorphism with open dense image.
It follows that the $G$-equivariant measure $\mu_t$ can be expressed as 
$w_t(Y) dY d(gM)$ with $w_t: \Omega^+\to \R_{>0}$ a measurable 
weight function.  

\par Now we note that that the right hand side of (\ref{3.3}) 
can be computed to be 

\begin{eqnarray} \notag  \int_\Xi \cK_x^t(z) \oline {\cK_y^t(z)} \ d\mu_t(z)
&=& \int_\Xi k_{2t}^\sim (a^{-1}z)  \oline{ k_{2t}^\sim (z)}\ d\mu_t(z) \\  
\notag &=& \int_G \int_{\Omega^+} 
k_{2t}^\sim (a^{-1}g\exp(iY)\cdot x_o)  \oline{ k_{2t}^\sim (g\exp(iY)\cdot x_o)} \, w_t(Y) \, dY\,  dg \\
\label{ab}&=& \int_G \int_{\Omega^+} 
k_{2t}^\sim (a^{-1}g\exp(iY)\cdot x_o)  k_{2t}^\sim (g\exp(-iY)\cdot x_o) \, w_t(Y) \, dY\,  dg. 
\end{eqnarray}

We perform an intermediate step. 
\begin{lemma}\label{aa} For all $a\in A$ and $Y\in \Omega$ the following 
identity holds:

\begin{align} \label{A}&  \int_G k_{2t}^\sim (a^{-1}g\exp(iY)\cdot x_o)  k_{2t}^\sim (g\exp(-iY)\cdot x_o) \ dg\\ 
\notag & \qquad = \int_{i\fa^*} 
e^{-4t(|\lambda|^2 +|\rho|^2)} \varphi_\lambda (a) \varphi_\lambda(\exp(2iY)\cdot x_o)
\ \frac{d\lambda}{ |\bc(\lambda)|^2} .
\end{align}
\end{lemma}

\begin{proof} Both sides of (\ref{A}) depend analytically on $Y$ and hence it 
is sufficient to prove the identity for $Y\in \frac{1}{2}\Omega$. 
Change of variables transforms the left hand side of 
(\ref{A}) 
into 
\begin{equation} \label{AB}
\int_G k_{2t}^\sim (a^{-1}g\exp(2iY)\cdot x_o)  k_{2t}(g\cdot x_o) \ dg
\end{equation}
and we start manipulating (\ref{AB}):

\begin{align*} & \int_G k_{2t}^\sim (a^{-1}g\exp(2iY)\cdot x_o)  k_{2t}(g\cdot x_o) \ dg\\ 
& \qquad  =
 \int_{i\fa^*} \int_G e^{-2t(|\lambda|^2 +|\rho|^2)}
\varphi_\lambda (a^{-1}g\exp(2iY)\cdot x_o) k_{2t}(g) \ dg \ \frac{d\lambda}{ |\bc(\lambda)|^2} \\ 
& \qquad =  \int_{i\fa^*} \int_G \int_K 
e^{-2t(|\lambda|^2 +|\rho|^2)} \varphi_\lambda (a^{-1}kg\exp(2iY)\cdot x_o) k_{2t}(g)\ dk 
 \ dg \ \frac{d\lambda}{|\bc(\lambda)|^2} \\ 
& \qquad =   \int_{i\fa^*} \int_G 
e^{-2t(|\lambda|^2 +|\rho|^2)} \varphi_\lambda (a^{-1}) \varphi_\lambda(g\exp(2iY)\cdot x_o) 
k_{2t}(g) \ dg \ \frac{d\lambda}{|\bc(\lambda)|^2}\\ 
& \qquad =   \int_{i\fa^*} 
e^{-4t(|\lambda|^2 +|\rho|^2)} \varphi_{-\lambda} (a) \varphi_{-\lambda}(\exp(2iY)\cdot x_o) 
\ \frac{d\lambda}{|\bc(\lambda)|^2} \\
& \qquad = \int_{i\fa^*} 
e^{-4t(|\lambda|^2 +|\rho|^2)} \varphi_\lambda (a) \varphi_\lambda(\exp(2iY)\cdot x_o)
\ \frac{d\lambda}{ |\bc(\lambda)|^2}\, . \end{align*}
\end{proof}

We combine Lemma \ref{aa}, (\ref{ab}) and (\ref{3.4}) to arrive at 
the curved analogue of the abelian identity (\ref{e2}). 

\begin{lemma} The existence of a $G$-invariant weight function for 
$\cF_t(\Xi)$ implies 

\begin{equation} \label{e3} k_{2t}(a)=\int_{\Omega^+} \int_{i\fa^*} 
e^{-4t(|\lambda|^2 +|\rho|^2)} \varphi_\lambda (a) 
\varphi_\lambda(\exp(2iY)\cdot x_o) w_t(Y) 
{\frac {d\lambda} {|\bc(\lambda)|^2}}\ dY  \quad (a\in A)\, .\end{equation} 
\end{lemma}

Finally, we derive a contradiction from (\ref{e3}). 
We replace $k_{2t}$ on the left hand side 
of (\ref{e3}) by its spectral resolution (\ref{eq=heat}), use uniqueness of the Fourier 
transform and conclude that (\ref{e3}) is equivalent to 

\begin{equation} \label{e4} e^{-2t(|\lambda|^2 +|\rho|^2)} =
  \int_{\Omega^+} e^{-4t(|\lambda|^2 +|\rho|^2)} 
\varphi_\lambda(\exp(2iY)\cdot x_o) w_t(Y) \, dY  \quad (\lambda\in i\fa^*)\,, \end{equation} 
or, in other words, 

\begin{equation} \label{e5} e^{2t(|\lambda|^2 +|\rho|^2)} =\int_{\Omega^+} 
\varphi_\lambda(\exp(2iY)\cdot x_o) w_t(Y) \, dY  \quad (\lambda\in i\fa^*)\, .\end{equation} 

At this point estimates for $\varphi_\lambda(\exp(2iY))$ 
slightly better than  
(\ref{eq=esa}) would allow us to obtain a contradiction. 
For complex groups, however, there is an alternative approach using Harish-Chandra's formula

\begin{equation} \varphi_\lambda(a)= \frac{{\bc}(\lambda)}{\delta(a)} \sum_{w\in W}\epsilon (w) e^{\lambda(w\log a)} \qquad (a\in A)\end{equation} 
with the functions $\delta(a)= (-1)^{|\Sigma^+|/2}\prod_{\alpha\in\Sigma^+}[e^ {\alpha(\log a)} - e^ {-\alpha(\log a)}] = (-1)^{|\Sigma^+|/2}\sum_{w\in W}
\epsilon (w) e^{\rho (w\log a)}$ and ${\bc}(\lambda)=\frac{\prod_{\alpha\in \Sigma^+} \la \rho, \alpha\ra} 
{\prod_{\alpha\in \Sigma^+} \la \lambda, \alpha\ra}$. 
Hence for $G$ complex, and setting $W_t(Y):= 
\delta (\exp(2iY))^{-1} w_t(Y)$ we derive 
from (\ref{e5}) the identity 

\begin{equation} e^{2t(|\lambda|^2 +|\rho|^2)} =
\int_{\Omega^+}{\bc}(\lambda)\sum_{w\in W} 
\epsilon (w)e^{\lambda(i2wY)} W_t(Y) \, dY  \quad (\lambda\in i\fa^*).
\end{equation} 

Since $\varphi_\lambda(w\cdot a)=\varphi_\lambda(a)$, we have $W_t(w \cdot Y)= \epsilon (w) W_t(Y)$. Then using the standard unfolding argument we get 

\begin{equation} e^{2t(|\lambda|^2 +|\rho|^2)} =\int_{\Omega} {\bc}(\lambda)
e^{\lambda(i2Y)} W_t(Y) \, dY  \quad (\lambda\in i\fa^*).\end{equation} 

From the definition of $\Omega$ it is clear that $\delta(a)$ is non-zero at the boundary of $\Omega$, so introduces no singularity to the integrand there. Whereas the smoothness of the spherical function shows there is no singularity of the integrand caused by the vanishing on $\delta(a)$ on root planes. Using the discussion on p. 329 of \cite{Wa}, we see that behaviour near the root planes contributes at most an additional factor of polynomial growth in $\lambda$; on the other hand, the integral away from the root planes can be estimated by exponential growth in $\lambda$ as $\Omega$ is compact. This contradicts that $e^{2t(|\lambda|^2 +|\rho|^2)}$ has exponential growth in $|\lambda|^2$ and completes the argument for complex groups. The general case should follow similarly using Gangolli's expansion for $\delta (a)\varphi_\lambda(a)$ and his estimates.

\section{The crown as a $G$-domain of holomorphy for a spherical function}

This section too is included to resolve a disagreement between the authors and the referee. The material is extracted  
from unpublished notes of the first and last named authors. 

\begin{theorem}\label{lasto-one} The crown is the maximal $G$-invariant 
domain in $X_\C$ to which a spherical function 
$\varphi_\lambda$ with $\lambda\in i\fa^*$ 
extends holomorphically. 
\end{theorem}

\begin{remark} We stress the importance of $G$-invariance in the 
statement of the theorem. Without this assumption 
the statement becomes false, see \cite{KS}, Th. 4.2. 
In order to obtain a feeling 
for this situation consider the function $f(z)=\sqrt {1-z}$ on the 
complex plane $\C$. Now the maximal ${\mathbb S}^1$-invariant domain 
of definition for $f$ is the unit disk $D=\{ z\in \C: |z|<1\}$ and 
there the function $f$ is bounded. Similarily, all non-trivial positive 
definite spherical functions are bounded on $\Xi$ by \cite{K} and the 
theorem asserts that they do not extend to a bigger 
$G$-invariant domain. 
\end{remark}

The theorem can be reduced to the basic case 
of $G=\Sl(2,\R)$ which will be presented first. 

\subsection{The basic case}

For this section we shall assume that $X=G/K=\Sl(2,\R)/\SO(2, \R)$. 
We follow the custom 
and make the standard choices 
$$A=\left\{ \begin{pmatrix} t & 0 \\ 0 & \frac{1}{t}\end{pmatrix}: 
t>0\right \} \quad \hbox{and}\quad 
N=\left \{ \begin{pmatrix} 1 & x \\ 0 & 1\end{pmatrix}: x\in \R\right \}\, .$$
Then 
$$\Omega=\left \{ \begin{pmatrix} x & 0 \\ 0 & - x \end{pmatrix}\in \fa: 
|x|<\frac{\pi }{4}\right\}=(-\pi/ 4, \pi/ 4)$$
Define a $K_\C$-invariant holomorphic function on 
$X_\C=\Sl(2,\C)/ \SO(2,\C)$ by  

$$P: X_\C \to \C, \ \ z\cdot x_o\mapsto \tr (z z^t)=
a^2 +b^2 +c^2 +d^2 \qquad (z=\begin{pmatrix} 
a & b\\  c& d\end{pmatrix}).$$
Of course, $P$ is nothing other than the elementary 
spherical function. 
\par For a convex $\cW$-invariant open set $\omega\subseteq 2\Omega$ 
one defines in \cite{K} the smallest $K_\C$-domain $\hat X_{\C, \omega}$
in  $X_\C$ which contains $A\exp(i\omega)\cdot x_o$. 
Specifically, we draw attention 
to the $K_\C$-domains
\begin{align}\hat X_{\C, \Omega}&=
\{ z\in X_\C : \RE  P(z)>0\}\, \\ 
\hat X_{\C, 2\Omega}&=
\{ z\in X_\C : P(z)\in \C \bs]-\infty, -2]\}\, .\end{align}

We recall the significance of $\hat X_{\C, 2\Omega}$ 
as the largest $K_\C$-domain to which spherical functions
extend, \cite{KS-II} Th. 2.4. 
Next, as $\Xi=G\exp(i\Omega)\cdot x_o$, it is immediate that 
\begin{equation}\label{K-inc} \Xi\subset \hat X_{\C,\Omega}
\end{equation}
and we note that (\ref{K-inc}) holds in full generality \cite{K}. 

\par The open interval $\Omega=(-\pi/4, \pi/ 4)$ is best possible for 
even the weaker inclusion 
\begin{equation}\label{K-inc2} \Xi\subset \hat X_{\C,2\Omega}
\end{equation}
to hold. The precise statement is as follows. 

\begin{lemma} Let $G=\Sl(2,\R)$. Then for 
$Y\in 2\Omega\bs \oline \Omega$, 
\begin{equation} 
G\exp(iY)\cdot x_0\nsubseteq \hat X_{\C, 2\Omega}.\end{equation}
More precisely, there exists a curve $\gamma(s)$, $s\in [0,1]$,  in 
$G$ such that the assignment  
$$s \mapsto  \sigma(s)=P(\gamma(s)\exp(iY)\cdot x_o)$$ 
is strictly decreasing with values in $[-2,2]$ such that 
$\sigma(0)=P(x_o)=2$ and 
$\sigma (1)=-2$.   
\end{lemma}

\begin{proof} Let $g=\begin{pmatrix} a &b\\  c& d\end{pmatrix}\in G$ and 
$z=\begin{pmatrix} e^{i\phi}& 0\\ 
0& e^{-i\phi} \end{pmatrix}\in \exp (2i\Omega)\setminus 
\exp(i\oline \Omega)$. 
This means $a,b,c,d\in \R$ with $ad-bc=1$ and 
$\frac{\pi}{ 4}<|\phi|<\frac{\pi}{ 2}$ for $\phi\in \R$.  
Thus 
\begin{eqnarray*} P(gz\cdot x_o)&=&P\begin{pmatrix} 
ae^{i\phi} & be^{-i\phi}\\ 
c e^{i\phi} &de^{-i\phi} \end{pmatrix}=
a^2e^{2i\phi} +b^2e^{-2i\phi} +c^2 e^{2i\phi} +d^2e^{-2i\phi}\cr 
&=&\cos (2\phi) (a^2+b^2 +c^2 +d^2) +i \sin 2\phi
(a^2-b^2 +c^2 -d^2)
\end{eqnarray*}

Using that $G=KAN$ and that $P$ is left $K$-invariant, we may actually 
assume that $g\in AN$, i.e. 
$$g=\begin{pmatrix} a& b\\  0 & \frac{1}{a}\end{pmatrix}$$
for some $a>0$ and $b\in \R$. Then 

$$P(gz\cdot x_o)=\cos (2\phi) (a^2+\frac{1}{a^2} +b^2) +i \sin 2\phi
(a^2-\frac{1}{a^2}-b^2 ).$$
We now show that $P(gz\cdot x_o)=-2$ has a solution for fixed  
$\frac{\pi}{ 4}<|\phi| <\frac{\pi}{ 2}$. 
This is because 
$P(gz\cdot x_o)=-2$ forces 
$\IM P(gz\cdot x_o)=0$ and so $b^2=a^2-\frac{1}{ a^2}$. 
Thus 
$$P(gz\cdot x_o)=2 a^2 \cos (2\phi) =-2.$$ 
Thus if we choose $a=\frac{1}{\sqrt{-\cos 2\phi}}$ we obtain a solution. 
The desired curve $\gamma(s)$ is now given by 
$$\gamma(s)=\begin{pmatrix} a(s) &b(s)\\  0& \frac{1}{a(s)}\end{pmatrix}$$
with $a(s)= \frac{1}{\sqrt{-\cos 2\phi}} 
(\sqrt{-\cos 2\phi} + s(1 -\sqrt{-\cos 2\phi})) $ and 
$b(s)=\sqrt{a(s)^2-\frac{1}{a(s)^2}}$. 
\end{proof}

\begin{theorem}\label{pope}  Let $G=\Sl(2,\R)$. Then 
the crown $\Xi$ is a maximal $G$-domain to which a 
spherical function  $\varphi_\lambda$ with  $\lambda\in i\fa^*$ 
extends holomorphically. 
\end{theorem}

\begin{proof}  Fix $\lambda\in i\fa^*$. 
We consider the spherical function $\varphi_\lambda$
on its  maximal $K_\C$-domain $\hat X_{\C, 2\Omega}$ of definition. 
Thus for each $\varphi_\lambda$ there exists 
a holomorphic function $\Phi_\lambda$ on $\C \bs (-\infty, 2]=
P(\hat X_{\C,2\Omega})$ 
such that 
\begin{equation} \label{mumu}\varphi_\lambda(z)=\Phi_\lambda(P(z)) \qquad 
(z\in \hat X_{\C, 2\Omega}).\end{equation}
Let $Y\in 2\Omega\bs \oline \Omega$. Let 
$\gamma\subset G$ and $\sigma\subset [-2,2]$ be curves as in the previous
lemma.
\par Note that $\gamma(s)\exp(iY)\cdot x_o\subset G$
for all $s\in [0,1)$. Hence (\ref{mumu}) gives 
$$\varphi_\lambda(\gamma(s)\exp(iY)\cdot x_o)=\Phi_\lambda(\sigma(s))
\qquad (s\in [0,1)\, .$$
Now recall that 
$s\mapsto \Phi_\lambda(\sigma(s))$ is positive (cf. \cite{KS}, Th. 4.2) 
and tends to infinity for $s\nearrow 1$ (cf.\ \cite{KS-II}, Th. 2.4).  
Thus the assertion of the theorem  will be proved if we can show
that a $G$-invariant domain $\Xi'\subset X_\C$ properly containing 
$\Xi$ contains points in $\exp(i2\Omega\bs \oline \Omega)$. 
But this follows from the fact that each $G$-orbit 
in the boundary $\partial \Xi$ contains the point 
$\begin{pmatrix}e^{i\pi/4} & 0 \\ 0 & e^{-i\pi/4} \end{pmatrix}\cdot x_o
\in \exp(i\partial \Omega)\cdot x_o$ in its closure
(cf. \cite{KS-II}, Lemma 2.3(ii)). \end{proof}

\subsection{Proof of Theorem \ref{lasto-one}}
Now that we understand $G=\Sl(2,\R)$, the general 
case will follow easily. 
We recall some facts on the boundary of $\Xi$.
First, each boundary $G$-orbit contains a point 
of $\exp(i\partial\Omega)\cdot x_o$ in its closure
(\cite{KS-II}, Lemma 2.3 (ii)). Next, to each $Y\in \partial \Omega$
we can associate an $\Sl(2,\R)$-crown  $\Xi_{\Sl_2}$
(cf. \cite{KS-II}, Th. 2.4
with formula (2.2) in its proof) which embeds
into $\Xi$ in such a way that 
$$\partial \Xi_{\Sl_2} 
\ni \begin{pmatrix}e^{i\pi/4} & 0 \\ 0 & e^{-i\pi/4} \end{pmatrix}\cdot x_o \mapsto \exp(iY)\cdot x_o\in \partial \Xi.$$
We restrict the spherical function $\phi_\lambda$ to $\Xi_{\Sl_2}$, 
see \cite{KS}, Prop. 4.5, and apply Theorem \ref{pope}.

\end{document}